\documentclass[aps,12pt,onecolumn,final,floats,bibnotes,showkeys]{revtex4}
\usepackage{amsfonts}
\usepackage{amsmath}
\usepackage{amssymb}
\usepackage{graphicx}

\begin{document}

\title{One and two side generalisations of the log-Normal distribution by
means of a new product definition}
\author{S\'{\i}lvio M. Duarte Queir\'{o}s $\circledS $}
\address{Unilever R\&D Port Sunlight\\
Quarry Road East, CH63 3JW, UK}
\email{sdqueiro@gmail.com}

\begin{abstract}
In this manuscript we introduce a generalisation of the log-Normal
distribution that is inspired by a modification of the Kaypten
multiplicative process using the $q$-product of Borges [Physica A \textbf{340%
}, 95 (2004)]. Depending on the value of q the distribution increases the tail for small (when $q<1$)
or large (when $q>1$) values of the variable
upon analysis. The usual log-Normal distribution is retrieved when $q=1$. The
main statistical features of this distribution are presented as well as a related random number
generators and tables of quantiles of the Kolmogorov-Smirnov.
Lastly, we illustrate the application of this distribution studying
the adjustment of a set of variables of biological and financial origin.
\end{abstract}


\date{The 28$^{th}$ August 2009}

\maketitle

\section{Introduction}

\label{intro}

The two-parameter log-Normal distribution, with probability density function,
\begin{equation}
p\left( x\right) =\frac{1}{\sqrt{2\pi }\sigma x}\exp \left[ -\frac{\left(
\ln x-\mu \right) ^{2}}{2\sigma ^{2}}\right] ,\qquad x>0,  \label{log-dist}
\end{equation}%
has played a major role in the statistical characterisation of many data sets for several decades
(empirical fitting) and has been an inspiration for theoretical studies as
well. The form of Eq.~(\ref{log-dist})
has been derived in several ways with particularly emphasis to the works of
Kapteyn \cite{kapteyn}, the Gibrat's law of proportionate effect~\cite{gibrat}, the Theory
of Breakage introduced by Kolmogorov~\cite{kolmogorov} or more recently in the theory of
chemical reactions~\cite{fa}. Concomitantly, Eq.~(\ref{log-dist}) has been
systematically modified to cope with different sets of data. Of those
generalisations the most famous of them is the truncated
log-normal distribution,%
\begin{equation}
p\left( x\right) =\frac{1}{\sqrt{2\pi }\sigma \left( x-\gamma \right) }\exp %
\left[ -\frac{\left( \ln \left( x-\gamma \right) -\mu \right) ^{2}}{2\sigma
^{2}}\right] ,\qquad 0<\gamma <x,  \label{log-truncated}
\end{equation}%
which has become a mathematical object of study in itself.

In this manuscript we introduce an alternative generalisation of the
log-normal distribution which we will term the $q$\emph{-log Normal}
distribution for historial reasons. This purported probability density
function emerges from changing the traditional algebra in the Kaypten
dynamics by a modified multiplication operation independently introduced by
Borges~\cite{borges} and Nivanen~\textit{et~al.}~\cite{nivanen}. This algebra has got
direct consequences on the emergence of asymptotic scale-free behaviour.
Specifically, in this manuscript we survey a new family of probability
density functions based on,%
\begin{equation}
p_{q}\left( x\right) =\frac{1}{\mathcal{Z}_{q}\,x^{q}}\exp \left[ -\frac{%
\left( \ln _{q}\,x-\mu \right) ^{2}}{2\,\sigma ^{2}}\right] ,\qquad \left(
x\geq 0\right) ,  \label{q-log}
\end{equation}%
where $ln_{q}(x)$ represents a generalisation of the logarithm of base $e$
with the normalisation,%
\begin{equation}
\mathcal{Z}_{q}=\left\{
\begin{array}{ccc}
\sqrt{\frac{\pi }{2}}\mathrm{erfc}\left[ -\frac{1}{\sqrt{2}\sigma }\left(
\frac{1}{1-q}+\mu \right) \right] \sigma  & if & q<1 \\
&  &  \\
\sqrt{\frac{\pi }{2}}\mathrm{erfc}\left[ \frac{1}{\sqrt{2}\sigma }\left(
\frac{1}{1-q}+\mu \right) \right] \sigma  & if & q>1.%
\end{array}%
\right. ,
\end{equation}
where $$\mathrm{erfc}(x) \equiv 2 \Phi(\sqrt{2} \, x)-1,$$ and
which in the limit $q\rightarrow 1^{\pm }$ exactly gives the traditional
log-Normal distribution. The aim of the present work is to introduce the functional
form of the distribution,
its dynamical origins and statistical features as well as applying it
to data of biological and financial origin. The manuscript is organised as
follows. In Sec. \ref{preliminar} we give a historical and mathematical
introduction of the underlying algebra; In Sec.~\ref{multiplicative} we reinterpret the Kaypten scenario for the emergence of
the log-normal, but using the $q$-algebra formalism which lead to the $q<1$,
$q>1$, and double $q$-log-Normal probability density function. In Secs.~\ref%
{exemplo} and \ref{randomnumber} we analyse their statistical properties and generate
random variables according to the distribution.
Finally, in Sec.~\ref{test} we introduce some real examples to which the new
distribution is shown to be a worthy candidate for medelling the data.

\section{Preliminaries: the $q$-product}

\label{preliminar}

The $q$-product, $\otimes _{q}$, has its origins in the endeavour to
extend the subject of statistical mechanics to systems exhibiting anomalous
behaviour when compared to systems described at Boltzmann-Gibbs equilibrium,
\textit{i.e.}, to deal with systems presenting long-lasting correlations,
ageing phenomena, non-exponential sensitivity to initial conditions, and
scale-invariance occupancy of the allowed phase space (for detailed
explanation of these concepts see~\cite{applications,abe}). The proposed
extension of statistical mechanics theory is grounded on the entropic
functional%
\begin{equation}
S_{q}\equiv \frac{1-\int \left[ p\left( x\right) \right] ^{q}\,dx}{q-1}%
,\qquad \left( q\in
\mathbb{R}
\right)  \label{ct-88}
\end{equation}%
(in its continuous and one-dimensional version) usually called Tsallis
entropy as well~\cite{ct-88}. This entropic form recovers the celebrated
Boltzmann-Gibbs-Shannon information measure,%
\begin{equation}
S=-\int p\left( x\right) \,\ln \left( x\right) \,dx,  \label{bg}
\end{equation}%
in the limit that the entropic parameter $q$ approaches $1$. The
interpretation of Eq. (\ref{ct-88}) as a $q$ generalisation of Eq.~(\ref{bg}%
) induced the introduction of analogue functions of the exponential and the
logarithm, namely, the $q$\emph{-exponential}
\begin{equation}
\exp _{q}\left( x\right) \equiv \left[ 1+\left( 1-q\right) \,x\right] ^{%
\frac{1}{1-q}},\qquad \left( x,q\in
\mathbb{R}
\right) ,
\end{equation}%
($\exp _{q}\left( x\right) =0$ if $1+(1-q)\,x\leq 0$) and its inverse the $q$%
\emph{-logarithm}~\cite{ct-quimica},
\begin{equation}
\ln _{q}\left( x\right) \equiv \frac{x^{1-q}-1}{1-q},\qquad \left( x>0,q\in
\mathbb{R}
\right) .
\end{equation}%
A functional form that generalises
the mathematical identity,
\begin{equation}
\exp \left[ \ln \,x+\ln \,y\right] =x\times y,\qquad \left( x,y>0\right) ,
\end{equation}%
for the $q$-product is,
\begin{equation}
x\otimes _{q}y\equiv \exp _{q}\left[ \ln _{q}\,x+\ln _{q}\,y\right] .
\label{q-product}
\end{equation}
For $q\rightarrow 1$, Eq.~(\ref{q-product}) recovers the usual
property,
\begin{equation*}
\ln \left( x\times y\right) =\ln \,x+\ln \,y
\end{equation*}%
($x,y>0$), with $x\times y\equiv x\otimes _{1}y$. Its inverse operation, the
$q$-division, $x\oslash _{q}y$, satisfies the following equality $\left(
x\otimes _{q}y\right) \oslash _{q}y=x$.

Bearing in mind that the $q$-exponential is a non-negative function, the $q$%
-product must be restricted to the values of $x$ and $y$ that respect the
condition,
\begin{equation}
\left\vert x\right\vert ^{1-q}+\left\vert y\right\vert ^{1-q}-1\geq0.
\label{cond-q-prod}
\end{equation}
Moreover, we can extend the domain of the $q$-product to negative values of $%
x$ and $y$ writing it as,
\begin{equation}
x\otimes_{q}y\equiv\mathrm{\ sign}\left( x\,y\right) \exp_{q}\left[ \ln
_{q}\,\left\vert x\right\vert +\ln_{q}\,\left\vert y\right\vert \right] .
\label{q-product-new}
\end{equation}
Regarding some key properties of the $q$-product we mention:

\begin{enumerate}
\item $x\otimes_{1}y=x\ y$;

\item $x\otimes_{q}y=y\otimes_{q}x$;

\item $\left( x\otimes_{q}y\right) \otimes_{q}z=x\otimes_{q}\left(
y\otimes_{q}z\right) =\left[ x^{1-q}+y^{1-q}-2\right] ^{\frac{1}{1-q}}$;

\item $\left( x\otimes_{q}1\right) =x$;

\item $\ln_{q}\left[ x\otimes_{q}y\right] \equiv\ln_{q}\,x+\ln_{q}\,y$;

\item $\ln_{q}\left( x\,y\right) =\ln_{q}\left( x\right) +\ln_{q}\left(
y\right) +\left( 1-q\right) \ln_{q}\left( x\right) \ln_{q}\left( y\right) $;

\item $\left( x\otimes_{q}y\right) ^{-1}=x^{-1}\otimes_{2-q}y^{-1}$;

\item $\left( x\otimes _{q}0\right) =\left\{
\begin{array}{ccc}
0 &  & \mathrm{if\ }\left( q\geq 1\ \mathrm{and\ }x\geq 0\right) \ \mathrm{%
or\ if\ }\left( q<1\ \mathrm{and\ }0\leq x\leq 1\right)  \\
&  &  \\
\left( x^{1-q}-1\right) ^{\frac{1}{1-q}} &  & \mathrm{otherwise}%
\end{array}%
\right. $

\end{enumerate}

For particular values of $q$, \textit{e.g.}, $q=1/2$, the $q$-product
provides nonnegative values at points for which the inequality $\left\vert
x\right\vert ^{1-q}+\left\vert y\right\vert ^{1-q}-1<0$ is verified.
According to the cut-off of the $q$-exponential, a value of zero for $%
x\otimes _{q}y$ is set down in these cases. Restraining our analysis of Eq.~(%
\ref{cond-q-prod}) to the sub-space $x,y>0$, we can observe that for $%
q\rightarrow -\infty $ the region $\left\{ 0\leq x\leq 1,0\leq y\leq
1\right\} $ is not defined. As the value of $q$ increases, the forbidden
region decreases its area, and when $q=0$, we have the limiting line given
by $x+y=1$, for which $x\otimes _{0}y=0$. Only for $q=1$, the entire set of $%
x$ and $y$ real values of has a defined value for the $q$-product. For $q>1$%
,\ the condition (\ref{cond-q-prod}) implies a region, $\left\vert
x\right\vert ^{1-q}+\left\vert y\right\vert ^{1-q}=1$ for which the $q$%
-product diverges. This undefined region augments its area as $q$ goes to
infinity. When $q=\infty $, the $q$-product is only defined in $\left\{
x\geq 0,0\leq y\leq 1\right\} \cup \left\{ 0\leq x\leq 1,y>1\right\} $.
Illustrative plots are presented in Fig.\ (1) of~\cite{part1}.

From the properties presented above we ascertain that the $q$-product has
got a neutral element and opposite and inverse elements under restrictions.
However, distributive property is not held and this fact thwarts the $q$%
-product of having commutative ring or field structures. Nevertheless, it
does not diminish the importance of this algebraic structure as other
algebras like the tropical algebra \cite{tropical}\ do not present all the
standard algebra properties and because the $q$-product represents a quite
rare case of a both-side non-distributive structure~\cite{green}.

Besides its inherent exquisiteness, this generalisation has found its own
field of applicability in the definition of the $q$-Fourier transform~\cite{sabir} which plays a key part in non-linear generalisations of the $q$%
-Central Limit Theorem~\cite{clt}, the definition of a modified characteristics methods
which allows the full analytical solution of the porous medium equations~\cite{sudq} and the structure of Pascal-Leibniz triangles~\cite{thierry}.

\section{Multiplicative processes as generators of distributions}

\label{multiplicative}

Multiplicative processes, particularly stochastic multiplicative processes,
have been the source of plentiful models applied in several fields of
science and knowledge. In this context, we can name the study of fluid
turbulence \cite{turbulence}, fractals \cite{feder}, finance \cite%
{mandelbrot}, linguistics \cite{murilinho}, etc. Specifically,
multiplicative processes play a very important role in the emergence of the
log-Normal distribution as a natural and ubiquitous distribution.
With regard to the dynamical origins of the log-Normal
distribution, we have mentioned in Sec. \ref{intro} the most celebrated
examples. Now, we shall give a brief account of the Kapteyn's process. To that, let us
consider a variable $\tilde{Z}$ obtained from a multiplicative random
process,%
\begin{equation}
\tilde{Z}=\prod\limits_{i=1}^{N}\tilde{\zeta}_{i},  \label{product}
\end{equation}%
where $\tilde{\zeta}_{i}$ are nonnegative microscopic variables associated
with a distribution $f^{\prime }\left( \tilde{\zeta}\right) $. If we
consider the following change of variables $Z\equiv \ln \tilde{Z}$, then we
have,%
\begin{equation*}
Z=\sum\limits_{i=1}^{N}\zeta _{i},
\end{equation*}%
with $\zeta \equiv \ln \tilde{\zeta}$. Assume now that $\zeta $ has a distribution $f\left( \zeta \right) $ with
mean $\mu $ and variance $\sigma ^{2}$. Then, $Z$ converges to the
Gaussian distribution in the limit of $N$ going to infinity as entailed by
the Central Limit Theorem~\cite{araujo}. Explicitly, considering that the
variables $\zeta $ are independently and identically distributed, the
Fourier Transform of $p\left( Z^{\prime }\right) $ is given by,

\begin{equation}
\mathcal{F}\left[ p\left( Z^{\prime }\right) \right] \left( k\right) =\left[
\int_{-\infty }^{+\infty }e^{i\,k\,\frac{\zeta }{N}}\,f\left( \zeta \right)
\,d\zeta \right] ^{N},  \label{fourier1}
\end{equation}%
where $Z^{\prime }=N^{-1}Z$. For all $N$, the integrand can be expanded as,%
\begin{equation}
\begin{array}{c}
\mathcal{F}\left[ p\left( Z^{\prime }\right) \right] \left( k\right) =\left[
\sum\limits_{n=0}^{\infty }\frac{\left( ik\right) ^{n}}{n!}\frac{%
\left\langle \zeta ^{n}\right\rangle }{N}\right] ^{N}, \\
\\
\mathcal{F}\left[ p\left( Z^{\prime }\right) \right] \left( k\right) =\exp
\left\{ N\ln \left[ 1+ik\frac{\left\langle \zeta \right\rangle }{N}-\frac{1}{%
2}k^{2}\frac{\left\langle \zeta ^{2}\right\rangle }{N^{2}}+O\left(
N^{-3}\right) \right] \right\} ,%
\end{array}%
\end{equation}%
where $\left\langle \zeta ^{n}\right\rangle$ represents the $n$th order raw moment of $\zeta $.
Expanding the logarithm,%
\begin{equation}
\mathcal{F}\left[ P\left( Z^{\prime }\right) \right] \left( k\right) \approx
\exp \left[ ik\mu -\frac{1}{2N}k^{2}\sigma ^{2}\right] .
\end{equation}%
Applying the inverse Fourier Transform, and reverting the $Z^{\prime }$
change of variables we finally obtain,%
\begin{equation}
p\left( Z\right) =\frac{1}{\sqrt{2\,\pi \,N}\sigma }\exp \left[ -\frac{%
\left( Z-N\,\mu \right) ^{2}}{2\,\sigma ^{2}\,N}\right] .
\end{equation}%
We can define the attracting distribution in terms of the original
multiplicative random process which yields the usual log-Normal distribution
\cite{log-normal-book},%
\begin{equation}
p\left( \bar{Z}\right) =\frac{1}{\sqrt{2\,\pi \,N}\sigma \,\bar{Z}}\exp %
\left[ -\frac{\left( \ln \bar{Z}-N\,\mu \right) ^{2}}{2\,\sigma ^{2}\,N}%
\right] .
\end{equation}

Although this distribution with two parameters, $\mu $ and $\sigma $, is
able to appropriately describe a large variety of data sets, there are cases
for which the log-Normal distribution fails statistical testing~\cite{log-normal-book}.
In some of these cases, such a failure has been overcome
by introducing different statistical distributions (e.g., Weibull
distributions~\cite{weibull, weibull1, weibull2}) or by changing the 2-parameter log-Normal
distribution into a 3-parameter log-Normal distribution~\cite{yuang, finney},%
\begin{equation}
p\left( x\right) =\frac{1}{\sqrt{2\,\pi }\sigma \,\left( x-\theta \right) }%
\exp \left[ -\frac{\left( \ln \left[ x-\theta \right] -\mu \right) ^{2}}{%
2\,\sigma ^{2}}\right] ,
\end{equation}%
which is very well characterised in the current scientific literature \cite%
{cohen}.

\subsection{One side generalisations}

\begin{figure}[tbh]

\includegraphics[width=0.45\columnwidth,angle=0]{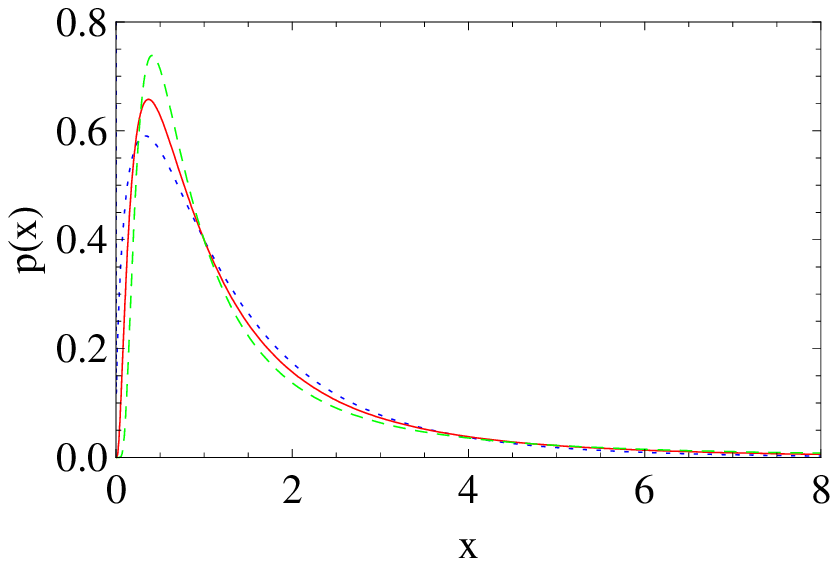} 
\includegraphics[width=0.45\columnwidth,angle=0]{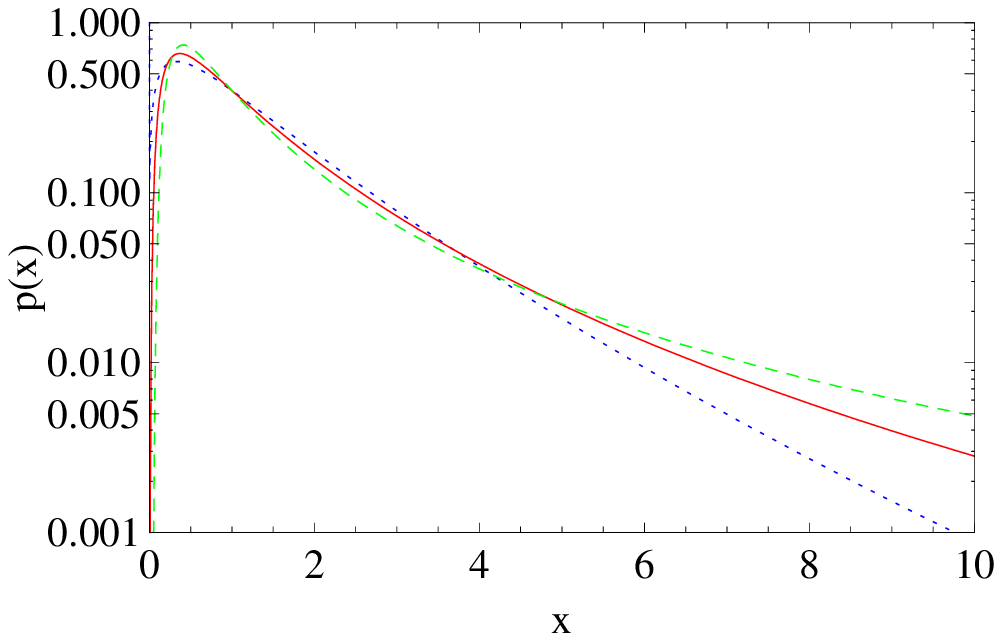} 
\includegraphics[width=0.45\columnwidth,angle=0]{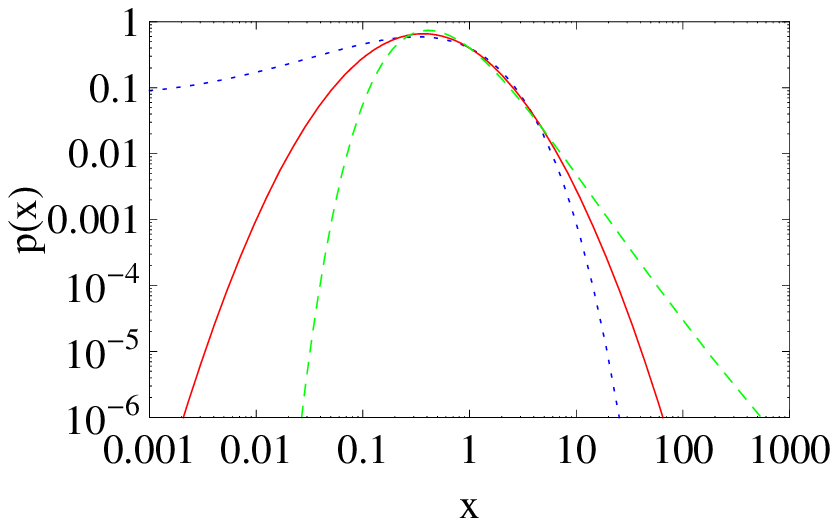}\\[0pt]

\caption{Plots of the Eq. (\protect\ref{qlog-normal}) \textit{vs} $x$ for $q=%
\frac {4}{5}$ (dotted line), $q=1$ (full line) and $q=\frac{5}{4}$ (dashed
line) in linear-linear scale (upper), log-linear (right), log-log (lower).}
\label{fig-pdf}
\end{figure}

Moving ahead, we now present our alternative procedure to generalise the distribution in
Eq.~(\ref{log-dist}). The motivation for this proposal comes from changing
the $N$ products in Eq. (\ref{product}) by $N$ $q$-products,%
\begin{equation}
\tilde{Z}=\underset{}{\prod\limits_{i=1}^{N}}^{(q)}\tilde{\zeta}_{i}\equiv
\tilde{\zeta}_{1}\otimes _{q}\tilde{\zeta}_{2}\otimes _{q}\ldots \otimes _{q}%
\tilde{\zeta}_{N}.  \label{p-product}
\end{equation}%
Applying the $q$-logarithm we have a sum of $N$ terms. If every term is
identically and independently distributed, then for variables $\zeta
_{i}=\ln _{q}$ $\tilde{\zeta}_{i}$ with finite variables we have a Gaussian
has stable distribution~\footnote{
Stable in the sense that if we consider the addition of two variables with that distribution the outcome of the convolution of the probability density functions is a probability density function with exactly the same functional form.}
, \textit{i.e.}, a Gaussian distribution in the $q$%
-logarithm variable. From this scenario we can obtain our $q$\emph{-log
Normal probability density function},%
\begin{equation}
p_{q}\left( x\right) \equiv \frac{1}{\mathcal{Z}_{q}\,x^{q}}\exp \left[ -%
\frac{\left( \ln _{q}\,x-\mu \right) ^{2}}{2\,\sigma ^{2}}\right] ,\qquad
\left( x\geq 0\right) ,  \label{qlog-normal}
\end{equation}%
with the normalisation,%
\begin{equation}
\mathcal{Z}_{q}\equiv \left\{
\begin{array}{ccc}
\sqrt{\frac{\pi }{2}}\mathrm{erfc}\left[ -\frac{1}{\sqrt{2}\sigma }\left(
\frac{1}{1-q}+\mu \right) \right] \sigma & if & q<1 \\
&  &  \\
\sqrt{\frac{\pi }{2}}\mathrm{erfc}\left[ \frac{1}{\sqrt{2}\sigma }\left(
\frac{1}{1-q}+\mu \right) \right] \sigma & if & q>1.%
\end{array}%
\right.
\end{equation}%
In the limit of $q$ equal to $1$, $\ln _{q\rightarrow 1}x=\ln x$ and $%
\mathcal{Z}_{q\rightarrow 1}=\sqrt{2\,\pi }\sigma $ and the usual log-Normal
is recovered. The cumulative distribution,
\begin{equation*}
\mathcal{P}\left( x\right) \equiv \int_{0}^{x}p\left( z\right) \,dz,
\end{equation*}%
is given by the following expressions,%
\begin{equation*}
\mathcal{P}_{q>1}\left( x\right) =\frac{1+\mathrm{erf}\left[ \frac{\ln
_{q}\left( x\right) -\mu }{\sqrt{2}\sigma }\right] }{1+\mathrm{erf}\left[
\frac{\frac{1}{q-1}-\mu }{\sqrt{2}\sigma }\right] },
\end{equation*}%
and,%
\begin{equation*}
\mathcal{P}_{q<1}\left( x\right) =\frac{\mathrm{erf}\left[ \frac{\ln
_{q}\left( x\right) -\mu }{\sqrt{2}\sigma }\right] -\mathrm{erf}\left[ -%
\frac{1}{\sqrt{2}\sigma }\left( \frac{1}{1-q}+\mu \right) \right] }{1+%
\mathrm{erfc}\left[ -\frac{1}{\sqrt{2}\sigma }\left( \frac{1}{1-q}+\mu
\right) \right] },
\end{equation*}%
Typical plots for cases with $q=\frac{4}{5}$, $q=1$, $q=\frac{5}{4}$
are depicted in Fig.~\ref{fig-pdf}. It can be seen that for $q$ greater than
one the likelihood of events round the peak as well as large values is
greater than that for the log-Normal case whereas the case $q<1$ favours events of
small value and the intermediate regime between the peak and the tail.

The raw statistical moments,%
\begin{equation}
\left\langle x^{n}\right\rangle \equiv \int_{0}^{\infty }x^{n}\,p\left(
x\right) \,dx,  \label{momentos}
\end{equation}%
can be analytically computed for $q<1$ giving \cite{gradshteyn},%
\begin{equation}
\left\langle x^{n}\right\rangle =\frac{\Gamma \left[ \nu \right] \exp \left[
-\frac{\gamma ^{2}}{8\,\beta }\right] D_{-\nu }\left[ \frac{\gamma }{\sqrt{%
2\,\beta }}\right] }{\sqrt{\beta ^{\nu }\,\pi }\,\sigma \left( 1-q\right)
\mathrm{erfc}\left[ -\frac{1}{\sqrt{2}\sigma }\left( \frac{1}{1-q}+\mu
\right) \right] },  \label{raw moment}
\end{equation}%
with%
\begin{equation}
\beta =\frac{1}{2\sigma ^{2}\left( 1-q\right) ^{2}};\quad \gamma =-\frac{%
1+\mu \,\left( 1-q\right) }{\left( 1-q\right) ^{2}\,\sigma ^{2}};\quad \nu
=1+\frac{n}{1-q},
\end{equation}%
where $D_{-a}\left[ z\right] $ is the \emph{parabolic cylinder function}
\cite{paraboliccylinderd}. Equation (\ref{raw moment}) allows us to write
the Fourier Transform or the generating function as,
\begin{equation}
\begin{array}{ccc}
\varphi \left( k\right) & = & \int p\left( x\right) e^{ikx}dx \\
&  &  \\
& = & \sum\limits_{n=0}^{\infty }\frac{\Gamma \left[ \nu _{n}\right] \exp %
\left[ -\frac{\gamma ^{2}}{8\,\beta }\right] D_{-\nu _{n}}\left[ \frac{%
\gamma }{\sqrt{2\,\beta }}\right] }{\sqrt{\beta ^{\nu _{n}}\,\pi }\,\sigma
\left( 1-q\right) \mathrm{erfc}\left[ -\frac{1}{\sqrt{2}\sigma }\left( \frac{%
1}{1-q}+\mu \right) \right] }(\mathrm{i}k)^{n}.%
\end{array}
\label{mom-gen}
\end{equation}

For $q>1$, the raw moments are given by an expression quite similar to the
Eq. (\ref{raw moment}) with the argument of the erfc replaced by
\begin{equation*}
\frac{1}{\sqrt{2}\sigma }\left( \frac{1}{1-q}+\mu \right) .
\end{equation*}%
However, the finiteness of the raw moments is not guaranteed for every $q>1$
for two reasons. First, according to the definition of $D_{-\nu }\left[ z\right] $, $%
\nu $ must be greater than $0$. Second, the core of the probability density
function (\ref{qlog-normal}),
\begin{equation*}
\exp \left[ -\frac{\left( \ln _{q}\,x-\mu \right) ^{2}}{2\,\sigma ^{2}}%
\right] ,
\end{equation*}%
does not vanish in the limit of $x$ going to infinity,%
\begin{equation}
\lim_{x\rightarrow \infty }\exp \left[ -\frac{\left( \ln _{q}\,x-\mu \right)
^{2}}{2\,\sigma ^{2}}\right] =\exp \left[ -\frac{\gamma ^{2}}{2}\right] .
\end{equation}%
This means that the limit $p\left( x\rightarrow \infty \right) =0$ is
introduced by the normalisation factor $x^{-q}$, which comes from redefining
the Normal distribution of variables,
\begin{equation}
y\equiv \ln _{q}x,
\end{equation}%
as the probability density function of variables $x$. Because of that, if
the moment exceeds the value of $q$, then the integral (\ref{momentos})
diverges. This has got severe repercussion in the adjustment procedures that
can be applied.

\subsection{Two side generalisation}

As visible from Fig.~(\ref{fig-pdf}), our generalisation modifies the tail
behaviour for small and large values of the variable depending on the value
of $q$ which describes the dependence between the variable that is
transferred into the $q$-product. It is well-known that many processes are
actually defined by a mixture of different laws of formation, some simpler
than others. Within this context, dual relations namely,%
\begin{eqnarray*}
q^{\prime } &=&2-q, \\
q^{\prime } &=&\frac{1}{2-q},
\end{eqnarray*}%
wherefrom property 7 of the $q$-product (see Sec. \ref{preliminar}) emerges,
are very inviting in the way that they represent the mapping of a certain rule
onto another which seems to be different at first but for which there is
actually a univocal transformation. Accordingly, we can imagine a scenario
in which variables follow two distinct paths either $q$-multiplying or $(2-q)
$-multiplying (corresponding to the inverse of the $q$-product) according to
some proportions $f$ and $f^{\prime }=1-f$. This proposal is in fact quite
plausible if we bear in mind few of the rife examples of mixing in dynamical
processes. From that, we establish the law,%
\begin{equation}
p_{q,2-q}\left( x\right) =f\,p_{q}\left( x\right) +\,f^{\prime
}\,p_{2-q}\left( x\right) ,  \label{two-branched}
\end{equation}%
for which we hold that $f=f^{\prime }=\frac{1}{2}$ is the most paradigmatic
case.

\subsection{Alternative interpretation}

The $q$-log-Normal distribution can introduce another clear advantage.
Namely, it provides us with an \emph{natural} and dynamical interpretation of the
truncated Normal distribution~\cite{johnsonkotz}. In other words, we can
look at the left(right) truncated Normal distribution,%
\begin{equation}
\mathcal{G}_{b}\left( y\right) =\sqrt{\frac{2}{\pi \,\sigma }}\text{\textrm{%
erfc}}\left[ \text{\textrm{sgn}}\left( b\right) \frac{1}{\sqrt{2}\sigma }%
\left( \mu -b\right) \right] ^{-1}\exp \left[ -\frac{\left( y-\mu \right)
^{2}}{2\,\sigma ^{2}}\right] ,
\end{equation}%
in which the truncation factor,
\begin{equation}
b=\frac{1}{q-1},  \label{truncation}
\end{equation}%
and $y=\ln _{q}x$ are intimately related to the value of $q$ which controls
the $q$-product part of the dynamical process. In this case the Fourier Transform can be analytically
determined. For left truncations we obtain,%
\begin{equation}
\mathcal{F}\left[ \mathcal{G}_{b}\right] \left( y\right) =\frac{1+\text{%
\textrm{erf\thinspace }}\left[ -\mathrm{i}\frac{k\sigma }{\sqrt{2}}-B\right]
}{\text{\textrm{erfc}}\left[ B\right] }\exp \left[ -\frac{1}{2}k\left( 2\,%
\mathrm{i}\,\mu +k\,\sigma ^{2}\right) \right] ,
\end{equation}%
where $\mathrm{erf}(x) \equiv \Phi(x)-1$. For right-truncations,%
\begin{equation}
\mathcal{F}\left[ \mathcal{G}_{b}\right] \left( y\right) =\frac{\text{%
\textrm{erfc\thinspace }}\left[ -\mathrm{i}\frac{k\sigma }{\sqrt{2}}-B\right]
}{\text{\textrm{erfc}}\left[ -B\right] }\exp \left[ -\frac{1}{2}k\left( 2\,%
\mathrm{i}\,\mu +k\,\sigma ^{2}\right) \right] ,
\end{equation}%
with%
\begin{equation}
B=\frac{b-\mu }{\sqrt{2}\sigma }.
\end{equation}%

\section{Examples of cascade generators}

\label{exemplo}

In this section, we discuss the upshot of two simple cases in which the
dynamical process described in the previous section is applied. We are going
to verify that the value of $q$ influences the nature of the attractor in
probability space.

\subsection{Compact distribution $[0,b]$}

\begin{figure}[tbh]

\includegraphics[width=0.75\columnwidth,angle=0]{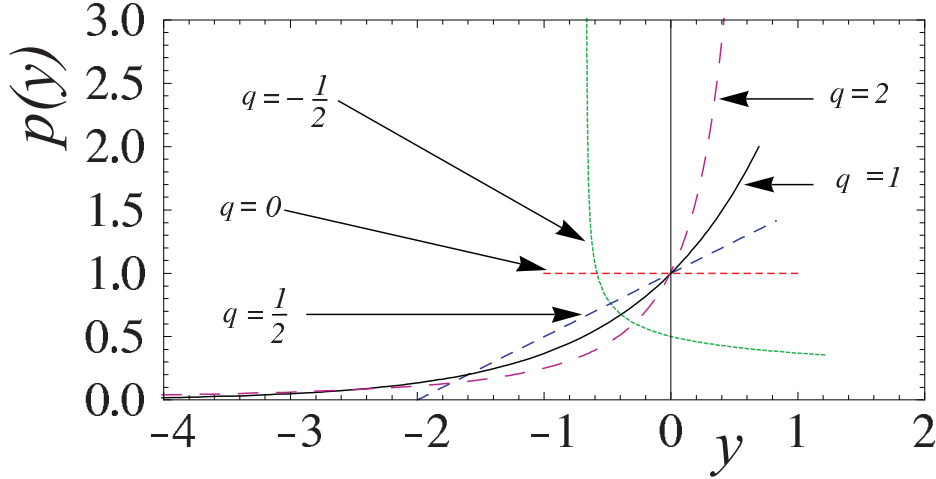}

\caption{Plots of the Eq. (\protect\ref{q-uniform}) \textit{vs} $y$ for $b=2$
and the values of $q$ presented in the text. }
\label{fig-flat}
\end{figure}

First, let us consider a compact distribution for identically and
independently distributed variables $x$ within the interval $0$ and $b$.
Following what we have described in the preceding section, we can transform
our generalised\ multiplicative process into a simple additive process of $%
y_{i}$ variables which are now distributed in conformity with the
distribution,%
\begin{equation}
p^{\prime}\left( y\right) =\frac{1}{b}\left[ 1+\left( 1-q\right) y\right] ^{%
\frac{q}{1-q}},  \label{q-uniform}
\end{equation}
with $y$ defined between $\frac{1}{q-1}$ and $\frac{b^{1-q}-1}{1-q}$ if $q<1$%
, whereas $y$ ranges over the interval between $-\infty$ and $\frac{b^{1-q}-1%
}{1-q}$ when $q>1$. Some curves for the special case $b=2$ are plotted in
Fig.~\ref{fig-flat}.

\begin{figure}[tbh]

\includegraphics[width=0.55\columnwidth,angle=0]{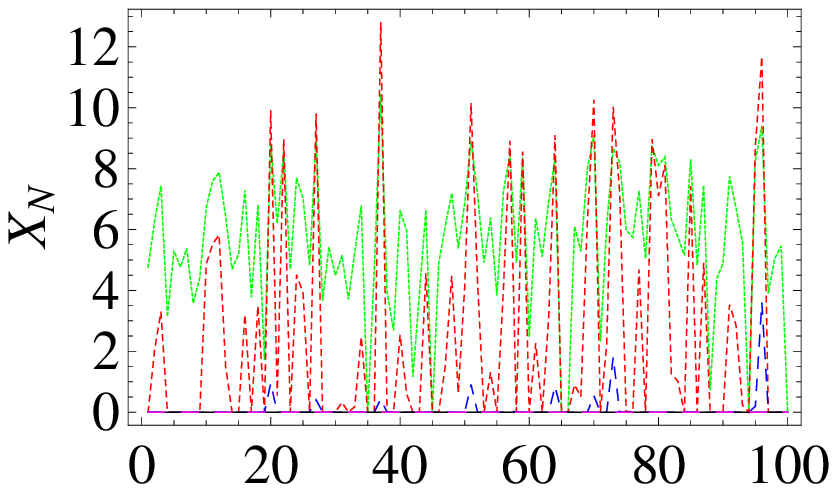} %
\includegraphics[width=0.55\columnwidth,angle=0]{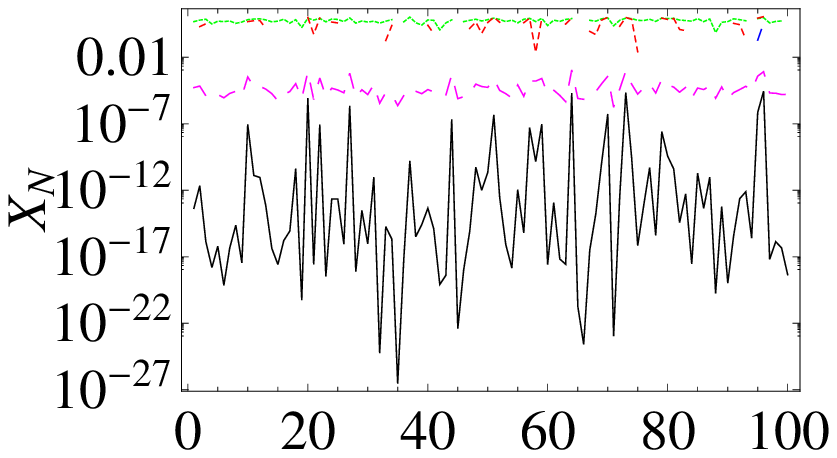}

\caption{Sets of random variables generated from the process (\protect\ref%
{p-product}) with $N=100$ and $q=-\frac{1}{2}$ (green), $0$ (red), $\frac{1}{%
2}$ (blue), $1$ (black), $\frac{5}{4}$ (magenta) in linear (upper panel) and
log scales (lower panel). The generating variable is uniformly distributed
within the interval $\left[ 0,1\right] $ as is the same for all of the cases
that we present. As visible, the value of $q$ deeply affects the values of $%
X_{N}=\tilde{Z}$. }
\label{fig-gen}
\end{figure}

If we look at the variance of this independent variable,%
\begin{equation}
\sigma _{y}^{2}=\left\langle y^{2}\right\rangle -\left\langle \mu
_{y}\right\rangle ^{2},
\end{equation}%
which is the moment whose finiteness plays the leading role in the Central
Limit Theory, we verify that for $q>\frac{3}{2}$, we obtain a divergent
value,%
\begin{equation}
\sigma _{y}^{2}=\frac{b^{2-2\,q}}{\left( 3-2\,q\right) \left( 2-q\right) ^{2}%
}.
\end{equation}%
Hence, if $q<\frac{3}{2}$, we can apply the Lyapunov's central Limit theorem
and our attractor in the probability space is the Gaussian distribution. On
the other hand, if $q>\frac{3}{2}$, the L\'{e}vy-Gnedenko's version of the
central limit theorem \cite{levy}\ asserts that the attracting distribution
is a L\'{e}vy distribution with a tail exponent,
\begin{equation}
\alpha =\frac{1}{q-1}.
\end{equation}%
Furthermore, it is simple to verify that the interval $\left( \frac{3}{2}%
,\infty \right) $ of $q$ values\ maps onto the interval $\left( 0,2\right) $
of $\alpha $ values, which is precisely the interval of validity of the L%
\'{e}vy class of distributions that is defined by its Fourier transform,%
\begin{equation}
\mathcal{F}\left[ L_{\alpha }\right] \left( k\right) =\exp \left[
-a\,\left\vert k\right\vert ^{\alpha }\right] .  \label{eq-levy}
\end{equation}%
In Fig.~\ref{fig-gen} we depict some sets generated by this process for
different values of $q$.

\subsection{$q$-log Normal distribution}

In this example, we consider the case of generalised multiplicative
processes in which the variables follow a $q$-log Normal distribution. In
agreement with what we have referred to in Sec. \ref{multiplicative}, the
outcome strongly depends on the value of $q$. Consequently, in the
associated $x$ space, if we apply the generalised process to $N$ variables $%
y=\ln _{q}\,x$ ($x\in \left[ 0,\infty \right) $) which follows a
Gaussian-like functional form \footnote{%
Strictly speaking, we cannot use the term Gaussian distribution because it
is not defined in the interval $\left( -\infty ,\infty \right) $. The
limitations in the domain do affect the Fourier transform and thus the
result of the convolution of the probability density function.} with average
$\mu $ and finite standard deviation $\sigma $, \textit{i.e.}, $\forall
_{q<1}$ or $q>3$ in Eq.(\ref{qlog-normal}), the resulting distribution in
the limit of $N$ going to infinity corresponds to the probability density
function (\ref{qlog-normal}) with $\mu \rightarrow N\,\mu $ and $\sigma
^{2}\rightarrow N\,\sigma ^{2}$. In respect of the conditions of $q$ we have
just mentioned here above, the $q$-log normal can be seen as an asymptotic
attractor, a stable attractor for $q=1$, and an unstable distribution for
the remaining cases with the resulting attracting distribution being
computed by applying the convolution operation.

\section{Random number generation and testing}

\label{randomnumber}

The generation of random numbers is by itself a subject of study and
depending on distribution different kinds of strategies can be used which
start from the shrewd von Neumann-Buffon acceptance-rejection method~\cite{vonNeuman}
and go to more sophisticated techniques~\cite{Gentle}. Since we
aim to introduce a global portrait of the distribution we have not tried to
develop bespoken algorithms but applied the robust method of the Smirnov
transformation (or inverse transformation sampling). From a classic robust
generator of uniform numbers, $\left\{ z\right\} $, between $-1$ and $1$ and
considering the probability conservation when $z$ is transformed into $y$
where $y$ is associated with a truncated log Normal distribution with
parameters $\mu $ and $\sigma $ and $b$ given by Eq. (\ref{truncation}). For
$q<1$, i.e., for $y=\ln _{q}\,x$ between $b$ and $+\infty $ we must use,%
\begin{equation}
y=\mu +\sqrt{2}\sigma \text{\thinspace \textrm{erf}\thinspace }^{-1}\left[
\infty ,\frac{1}{2}\left( z-1\right) \text{\textrm{erfc}}\left[ B\right] %
\right] ,  \label{generetormaior}
\end{equation}%
whereas for $q>1$, i.e., for $y=\ln _{q}\,x$ between $-\infty $ and $b$%
\begin{equation}
y=\mu +\sqrt{2}\sigma \text{\thinspace \textrm{erf}\thinspace }^{-1}\left[ 0,%
\frac{1}{2}\left\{ z\left( \text{\textrm{erf\thinspace }}\left[ B\right]
+1\right) +\text{\textrm{erf\thinspace }}\left[ B\right] -1\right\} \right] .
\label{generatormenor}
\end{equation}%
From these formulae we have defined the Kolmogorov-Smirnov distance tables
that we present for typical cases $q=4/5$ and $q=5/4$ with $\mu =0$ and $%
\sigma =1$. For each case $10^{6}$ samples have been considered.

\begin{table}
\caption{Quantiles of the Kolmogorov-Smirnov statistics of a q-log-Normal
distribution with $q=4/5$, $\protect\mu =0$ and $\protect\sigma =1$.}
\label{tab1}
\centering
\fbox{%
\begin{tabular}{*{6}{c}}
\hline
$n$ $\diagdown $ $P$ & 0.80 & 0.85 & 0.90 & 0.95 & 0.99 \\ \hline\hline
\multicolumn{1}{l}{5} & \multicolumn{1}{l}{0.442} & \multicolumn{1}{l}{0.471}
& \multicolumn{1}{l}{0.504} & \multicolumn{1}{l}{0.558} & \multicolumn{1}{l}{
0.663} \\
\multicolumn{1}{l}{10} & \multicolumn{1}{l}{0.318} & \multicolumn{1}{l}{0.339
} & \multicolumn{1}{l}{0.362} & \multicolumn{1}{l}{0.404} &
\multicolumn{1}{l}{0.485} \\
\multicolumn{1}{l}{15} & \multicolumn{1}{l}{0.261} & \multicolumn{1}{l}{0.277
} & \multicolumn{1}{l}{0.296} & \multicolumn{1}{l}{0.333} &
\multicolumn{1}{l}{0.399} \\
\multicolumn{1}{l}{20} & \multicolumn{1}{l}{0.211} & \multicolumn{1}{l}{0.225
} & \multicolumn{1}{l}{0.245} & \multicolumn{1}{l}{0.274} &
\multicolumn{1}{l}{0.334} \\
\multicolumn{1}{l}{25} & \multicolumn{1}{l}{0.192} & \multicolumn{1}{l}{0.205
} & \multicolumn{1}{l}{0.222} & \multicolumn{1}{l}{0.249} &
\multicolumn{1}{l}{0.302} \\
\multicolumn{1}{l}{30} & \multicolumn{1}{l}{0.176} & \multicolumn{1}{l}{0.188
} & \multicolumn{1}{l}{0.204} & \multicolumn{1}{l}{0.228} &
\multicolumn{1}{l}{0.278} \\
\multicolumn{1}{l}{35} & \multicolumn{1}{l}{0.164} & \multicolumn{1}{l}{0.175
} & \multicolumn{1}{l}{0.190} & \multicolumn{1}{l}{0.212} &
\multicolumn{1}{l}{0.257} \\
\multicolumn{1}{l}{40} & \multicolumn{1}{l}{0.154} & \multicolumn{1}{l}{0.165
} & \multicolumn{1}{l}{0.178} & \multicolumn{1}{l}{0.200} &
\multicolumn{1}{l}{0.242} \\
\multicolumn{1}{l}{45} & \multicolumn{1}{l}{0.146} & \multicolumn{1}{l}{0.156
} & \multicolumn{1}{l}{0.169} & \multicolumn{1}{l}{0.189} &
\multicolumn{1}{l}{0.230} \\
\multicolumn{1}{l}{50} & \multicolumn{1}{l}{0.140} & \multicolumn{1}{l}{0.149
} & \multicolumn{1}{l}{0.161} & \multicolumn{1}{l}{0.18} &
\multicolumn{1}{l}{0.219} \\
\multicolumn{1}{l}{60} & \multicolumn{1}{l}{0.128} & \multicolumn{1}{l}{0.139
} & \multicolumn{1}{l}{0.148} & \multicolumn{1}{l}{0.165} &
\multicolumn{1}{l}{0.199} \\
\multicolumn{1}{l}{70} & \multicolumn{1}{l}{0.119} & \multicolumn{1}{l}{0.127
} & \multicolumn{1}{l}{0.138} & \multicolumn{1}{l}{0.154} &
\multicolumn{1}{l}{0.186} \\
\multicolumn{1}{l}{80} & \multicolumn{1}{l}{0.112} & \multicolumn{1}{l}{0.121
} & \multicolumn{1}{l}{0.1301} & \multicolumn{1}{l}{0.144} &
\multicolumn{1}{l}{0.175} \\
\multicolumn{1}{l}{90} & \multicolumn{1}{l}{0.106} & \multicolumn{1}{l}{0.113
} & \multicolumn{1}{l}{0.122} & \multicolumn{1}{l}{0.135} &
\multicolumn{1}{l}{0.164} \\
\multicolumn{1}{l}{100} & \multicolumn{1}{l}{0.101} & \multicolumn{1}{l}{
0.108} & \multicolumn{1}{l}{0.115} & \multicolumn{1}{l}{0.129} &
\multicolumn{1}{l}{0.156} \\
$n>100$ & 1.02 $n^{-1/2}$ & 1.17 $n^{-1/2}$ & 1.30 $n^{-1/2}$ & 1.42 $%
n^{-1/2}$ & 1.56 $n^{-1/2}$ \\ \hline
\end{tabular}}%

\end{table}

\begin{table}
\caption{Quantiles of the Kolmogorov-Smirnov statistics of a q-log-Normal
distribution with $q=5/4$, $\protect\mu=0$ and $\protect\sigma = 1$.}
\label{tab2}
\centering
\fbox{%
\begin{tabular}{*{6}{c}}
\hline
$n$ $\diagdown$ $P$ & 0.80 & 0.85 & 0.90 & 0.95 & 0.99 \\ \hline\hline
5 & 0.382 & 0.413 & 0.454 & 0.513 & 0.627 \\
10 & 0.286 & 0.307 & 0.334 & 0.377 & 0.461 \\
15 & 0.246 & 0.262 & 0.282 & 0.317 & 0.384 \\
20 & 0.204 & 0.218 & 0.237 & 0.277 & 0.327 \\
25 & 0.189 & 0.202 & 0.218 & 0.246 & 0.299 \\
30 & 0.174 & 0.186 & 0.202 & 0.225 & 0.276 \\
35 & 0.161 & 0.174 & 0.188 & 0.213 & 0.256 \\
40 & 0.155 & 0.165 & 0.178 & 0.204 & 0.242 \\
45 & 0.146 & 0.156 & 0.169 & 0.189 & 0.229 \\
50 & 0.137 & 0.148 & 0.162 & 0.183 & 0.217 \\
60 & 0.128 & 0.138 & 0.148 & 0.165 & 0.201 \\
70 & 0.118 & 0.127 & 0.138 & 0.154 & 0.186 \\
80 & 0.111 & 0.121 & 0.1301 & 0.143 & 0.175 \\
90 & 0.107 & 0.113 & 0.122 & 0.135 & 0.164 \\
100 & 0.099 & 0.107 & 0.115 & 0.128 & 0.155 \\
$n>100$ & 1.01 $n ^{-1/2}$ & 1.15 $n ^{-1/2}$ & 1.28 $n ^{-1/2}$ & 1.41 $n
^{-1/2}$ & 1.57 $n ^{-1/2}$ \\ \hline
\end{tabular}}%

\end{table}

\section{Examples of applicability}

\label{test}

In the following examples parameter estimation has been made using
traditional maximum log-likelihood methods. In spite of using Brent's method
for optimisation~\cite{brent} of the log-likelihood function, the following
set of equations can be solved if a differential method is preferred:

\begin{equation}
\left\{
\begin{array}{c}
\sum\limits_{i=i}^{n}\frac{d}{dq}\ln P\left( x_{i}\right) =0 \\
\sum\limits_{i=i}^{n}\frac{d}{d\mu }\ln P\left( x_{i}\right) =0 \\
\sum\limits_{i=i}^{n}\frac{d}{d\sigma }\ln P\left( x_{i}\right) =0%
\end{array}%
\right. .
\end{equation}%
The specific equations can be obtained after straightforward (and tedious)
calculus.

\subsection{Shadow prices in metabolic networks}

The representation of metabolic networks is often related to linear
programming approaches~\cite{palsson} for which there is a dual optimisation
procedure. In other words, the maximisation of the reaction fluxes of a
metabolic network with a given stoichiometry matrix has as its dual
solution the minimisation of a certain function defined by quantities traditionally called
\textit{shadow prices}, $\Pi $, which for this case correspond to the chemical
potencial~\cite{pbw}. In a previous study the shape of the distribution of the shadow
prices has been analysed. From the set of tested PDFs the log-normal has
proven to be the better description.

\begin{figure}[tbh]

\includegraphics[width=0.75\columnwidth,angle=0]{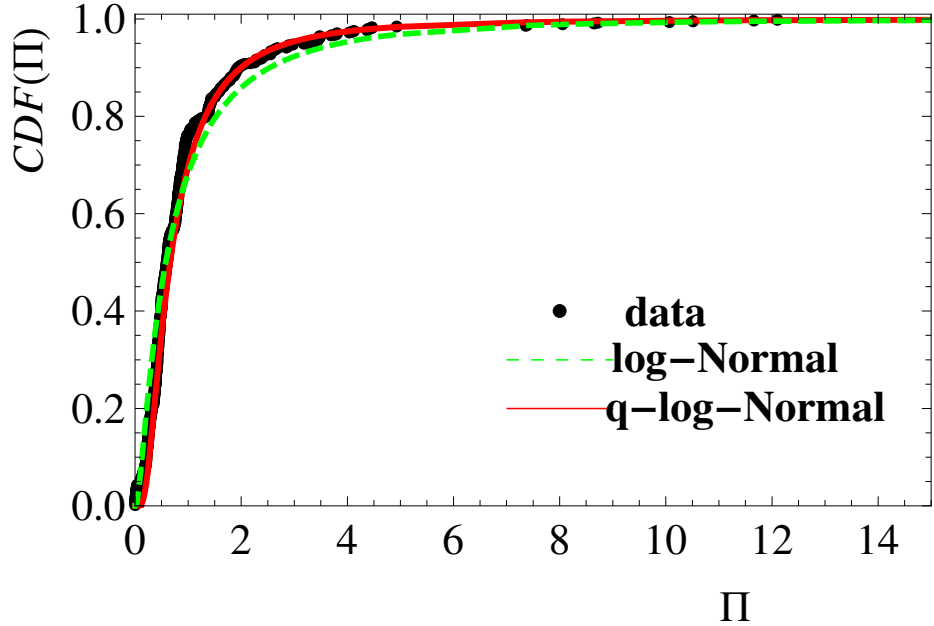}

\caption{Cumulative density function of the shadow prices vs shadow price of
the metabolic network of the E. coli (iJR 904) growing on a D-glucose
substrate. The symbols are obtained from the data and the lines the best
fits with the q-log-Normal distribution and log-Normal. The values of the
parameters and error are mentioned in the text.}
\label{Ecoli}
\end{figure}

Our first example is composed of shadow prices of the genome-scale model for
\textit{E. coli} (iJR 904) growing on a D-glucose substrate \cite{reed,kummel,reed2}. The number of
shadow prices is $649$. Minimisation of the log-likelihood function we
obtained $\mu =-0.432$, $\sigma =0.838$ and $q=1.21$ in comparison with $\mu
=-0.454$ and $\sigma =0.741$ for the log-Normal. The corresponding
Kolmogorov-Smirnov distances are $0.072$ and $0.123$, respectively, as we
depict in Fig.\ref{Ecoli}. Other qualitatively similar results, $q>1$, are
found for the shadow prices of models growing upon aerobic conditions.

A different kind of distribution was obtained when a metabolic network like
the \textit{M. barkeri} (iAF 692 model) evolving in a Hydrogen medium was
considered. In this case $517$ metabolites are taken into account. The
values of the best fit obtained were $\mu =-0.633$, $\sigma =1.24$ and $%
q=0.822$ in comparison with $\mu =-0.454$ and $\sigma =0.741$ for the
log-Normal.

\begin{figure}[tbh]

\includegraphics[width=0.75\columnwidth,angle=0]{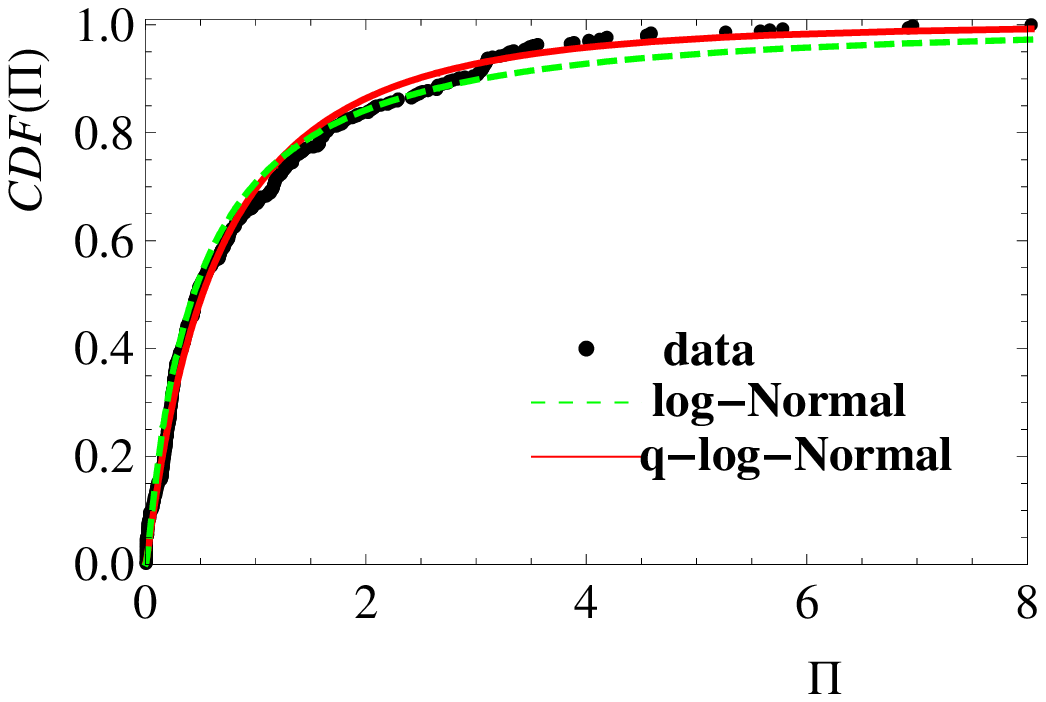}

\caption{Cumulative density function of the shadow prices vs shadow price of
the metabolic network of the M. barkeri (iAF 692 model) growing in Hydrogen
medium \cite{feist}. The symbols are obtained from the data and the lines the best fits
with the q-log-Normal distribution and log-Normal. The values of the
parameters and error are mentioned in the text.}
\label{MBarkeri}
\end{figure}

These parameters yield the following Kolmogorov-Smirnov distances $0.049$
and $0.080$ which represent a blunt improvement. Moreover, the introduction
of the extra parameter $q$ is completely justified when we calculate the
Akeike information criterion (AIC) \cite{aic},%
\begin{equation*}
AIC=2\,k+n\ln \left[ \frac{RSS}{n}\right] ,
\end{equation*}%
where $k$ is the number of parameters, $n$ is the number of metabolites of
the metabolic network and $RSS$ is the residual sum of squares. The values
of AIC per metabolite are $-6.288$ and $-5.607$ for the $q$-log-Normal and
the log-Normal in the case of the \textit{E. coli}, respectively. For the
\textit{M. barkeri} the values are $-7.474$ and $-6.682$. This is clearly
adduces that the $q$-log-Normal outperforms the log-Normal distribution
which had given the best results. From a biological perspective is even more
appealing that metabolic networks developed in a aerobic environment have
presented a value of $q>1$ and \ networks related to anaerobic environments
yield $q$ values smaller than 1. Whence, we can infer that $q$ value can be
possibly used as a signature of aerobic and anaerobic growing.

\subsection{Volatility in financial markets}

\begin{figure}[tbh]

\includegraphics[width=0.75\columnwidth,angle=0]{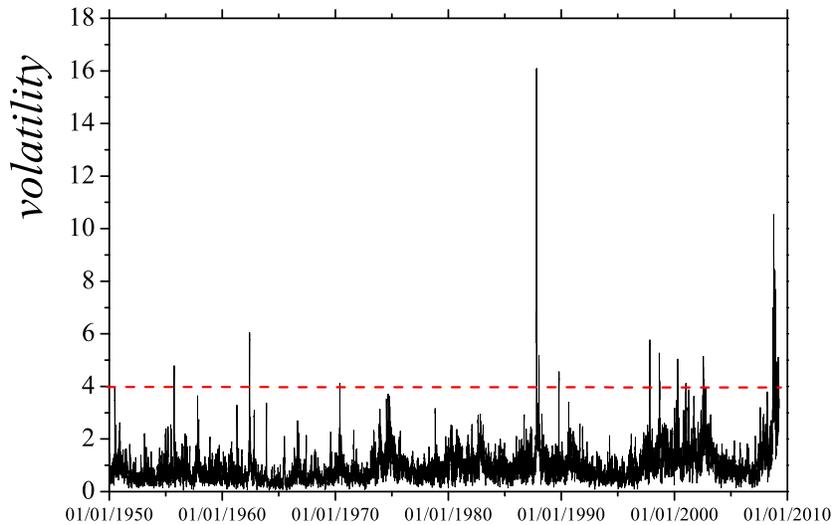}

\caption{Evolution of the five-day volatility of the SP500 index as defined
in the text after normalisation by its average value. Value above the dashed
red line can be considered extreme events.}
\label{volatilitySP}
\end{figure}

One of the keystone elements of mathematical finance is the volatility.
Despite appearing in every theory of financial markets the truth is that
volatility still lacks a precise definition~\cite{engle}. Nonetheless, it is
customarily associated with average of squared fluctuations,%
\begin{equation*}
r_{\Delta }\left( t\right) \equiv \ln S\left( t+\Delta \right) -\ln S\left(
t\right) ,
\end{equation*}%
of the (log-)price (or index) $S$ over some window $T$ ($\Delta $ is lag).
It is well-known for a long time that price fluctuations are nicely fitted by
the Student's $t$-distribution. An explanation for that relies on the local
Gaussianity of the price distributions but with a time dependent variance as
it has been hold by heteroskedastic processes~\cite{engle}. In that
sense, we can consider a variable $\mathcal{B}\left( t\right) =v\left(
t\right) ^{-1}$, where%
\begin{equation}
v\left( t\right) =\frac{1}{T}\sum_{i=1}^{T}r^{2}\left( t-i\right) .
\label{volatility}
\end{equation}%
Accordingly, the distribution of price fluctuations following a Bayesian
approach%
\begin{equation}
p\left( r\right) =\int P\left( \beta \right) \,p\left( r|\mathcal{B}\right)
\,d\mathcal{B}.  \label{superstatistics}
\end{equation}%
Assuming the Student's $t$-distribution hypothesis for $p\left( r\right) $ and
\begin{equation}
p\left( r|\mathcal{B}\right) =\sqrt{\frac{\mathcal{B}}{2\pi }}\exp [-%
\mathcal{B}\,r^{2}],  \label{localgauss}
\end{equation}%
the distribution of $\mathcal{B}$ must be,%
\begin{equation}
P_{\Gamma }\left( \mathcal{B}\right) =\frac{\theta ^{-1-\kappa }}{\Gamma %
\left[ 1+\kappa \right] }\mathcal{B}^{\kappa }\exp [-\frac{\mathcal{B}}{%
\theta }].  \label{p-gamma}
\end{equation}%
In this case, we have analysed the distribution of $\mathcal{B}$ according
to the definition given in Eq. (\ref{volatility}) using the daily
fluctuations of the SP500 index from the $3^{rd}$ January 1950 to the $%
3^{rd} $ April 2009 and considering 5-business days windows with data
obtained from {\tt http://finance.yahoo.com} (see Fig. \ref{volatilitySP}).

\begin{figure}[tbh]

\includegraphics[width=0.75\columnwidth,angle=0]{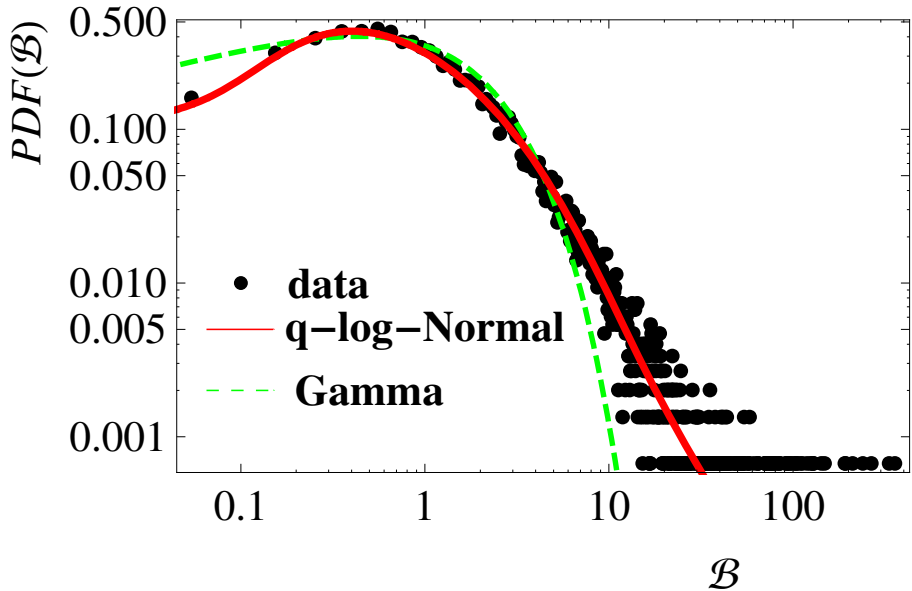}

\caption{Probability density function of the 5-day volatility vs $\mathcal{B}
$ . The symbols are obtained from the data and the lines the best fits with
the Gamma distribution and the double sided q-log-Normal. The values of the
parameters and error are mentioned in the text.}
\label{pdf-vol}
\end{figure}

Regarding $P(\mathcal{B})$, we have actually noticed that the distribution
is poorly described by Eq.~(\ref{p-gamma}). Conversely, we verified a very
good agreement with Eq.~(\ref{two-branched}) as we show in Fig.~\ref{pdf-vol}%
. The values obtained for $P_{\Gamma }\left( \mathcal{B}\right) $ are $%
\kappa =0.314$ and $\theta =1.41$ and for $p_{q,2-q}\left( \mathcal{B}%
\right) $ we have $\mu =0.391$, $\sigma =1.15$ and $q=1.22$. The values
of the Kolmogorov-Smirnov distances yield 0.0959 and 0.0126 which has passed
the statistical test for $\alpha =2\%$. A representation of the probability
density function adjustment is shown in Fig. \ref{pdf-vol}. Although not
shown here a log-Normal adjustment which yielded $\mu =0.379$ and $\sigma
=1.121$ and a Kolmogorov-Smirnov distance of $0.0177$ which is $40\%$
greater than the Kolmogorov-Smirnov distance of the $q$-Log-Normal. The
utilisation of two values for $q$ (although they relate one another) can be
understood if we accept tested hypotheses that the volatility runs over two
mechanisms (short and long scale) \cite{bouchaud}. Nevertheless, it is
worth mentioning that applying Eqs. (\ref{superstatistics})-(\ref%
{p-gamma}) with the values we have determined brings about a Student-$t$
distribution with $\nu =2.64$ that is in accordance with the value measured
for the tail exponent of that distribution (see e.g. \cite{smdqQF}).

\section{Final remarks}

In this manuscript we have introduced a new kind of generalisation of the
log-Normal distribution which stem from a modification of the multiplicative
cascade using a new type of algebra recently introduced in a physical
context. This modification of the $q$-product, as shown in other cases,
represents a way of describing a type of dependence between the variables.
Accordingly, the new distribution differs from the classical log-Normal by a
single parameter $q$ which favours the right-side tail for $q>1$, the
left-side tail if $q<1$ and recovers the traditional form when $q=1$. We
have made an extensive description of the distribution namely by defining
the moments, the Fourier Transform we have purported a generator of random
numbers as well which yield the distributions we mentioned here. Using these
random number generators we have depicted the construction of a $P$-value
table for the Kolmogorov-Smirnov distance when the $q$-log-Normal
distribution is assumed. Moreover, we have tested the distribution against
real data of biological and financial origin. Both results have shown its
usefulness and all the cases we have studied curiously present values of $q$
or $2-q$ close to $5/4$. Concerning future work we can mention the
modification of the two branched distribution to accommodate equal weights
for $f$ and $f^{\prime }$ as we have considered here, using different dual
relations for the $q$ parameters or parameters that are not related by any
dual relation as well. This last approach corresponds to accepting different
mixtures of dynamical or structural processes. It is obvious that such
modifications augment the number of the parameters which might be plainly
justified by usual statistical criteria.

\bigskip

{\small SMDQ acknowledges P.B. Warren for having provided the shadow prices data,
T. Cox for several comments on the subject matter and the critical reading of the manuscript
and M.A. Naeeni for preliminary discussions.
This work benefited from financial support from the Marie Curie Fellowship programme (European
Union).}


\begin{thebibliography}{}

\bibitem[{Abe {\it  et al.}}{, 2007}]{abe} (2007) \textit{Complexity, Metastability, and Nonextensivity: An International
Conference} (eds S. Abe, H. Herrmann, P. Quarati, A. Rapisarda, C. Tsallis), AIP Conf. Proc., \textbf{965}.

\bibitem[{Akeike}{, 1974}]{aic} Akeike H. (1974)
A new look at the statistical model identification. \textit{IEEE Trans. Aut. Control}, \textbf{19}, 716.

\bibitem[{Araujo and Guin\'{e}}{, 1980}]{araujo} Araujo A. and Guin\'{e} E. (1980)
\textit{The Central Limit Theorem for Real and Banach Valued Random Variables}. New York:John Wiley \& Sons.

\bibitem[{Beck et al.}{, 2005}]{bcs} Beck C., Cohen E.G.D., and Swinney H.L. (2005)
From time series to superstatistics. \textit{Phys. Rev. E}, \textbf{72}, 056133.

\bibitem[{Borges}{, 2004}]{borges} Borges E. P. (2004) A possible deformed algebra and
calculus inspired in nonextensive thermostatistics. \textit{Physica A}, \textbf{340}, 95.

\bibitem[{Bouchaud and Potters}{, 2000}]{bouchaud} Bouchaud J. P. and Potters M. (2000)
\textit{Theory of financial risk and derivative pricing: from statistical physics to risk management}. Cambridge: Cambridge University Press.

\bibitem[{Brent}{, 1973}]{brent} Brent R. P. (1973)
\textit{Algorithms for Minimization Without Derivatives}. Englewood Cliffs - NJ: Prentice \& Hall.

\bibitem[{Cohen}{, 1988}]{cohen} A.C. Cohen (1988) Three-parameter estimation in Lognormal Distributions: Theory and Applications. In
\textit{Lognormal Distributions: Theory and Applications} (eds E.L. Crow and K. Shimizu). New York: CRC Pess.

\bibitem[{Clay}{, 1992}]{tropical} Clay J. R. (1992) \textit{Nearrings: Genesis and Applications}. Oxford:Oxford University Press.

\bibitem[{Crow and Shimizu}{, 1988}]{log-normal-book} (1988) \textit{Lognormal Distributions: Theory and
Applications} (eds E.L. Crow and K. Shimizu). New York: CRC Pess.

\bibitem[{Duarte Queir\'{o}s}{, 2005}]{smdqQF} Duarte Queir\'{o}s S. M. (2005)
On non-Gaussianity and dependence in financial time series: a nonextensive approach. \textit{Quant. Finance} \textbf{5}, 475.

\bibitem[{Duarte Queir\'{o}s and Tsallis},{2007}]{part1} Duarte Queir\'{o}s S. M. and Tsallis C. (2007) Nonextensive statistical
mechanics and central limit theorems I - Convolution of independent random variables and the q-product. In
\textit{Complexity, Metastability, and Nonextensivity: An International Conference} (eds S. Abe, H. Herrmann,
P. Quarati, A. Rapisarda, C. Tsallis), AIP Conf. Proc. \textbf{965}, 8.

\bibitem[{Engle}{, 1995}]{engle} (1995) \textit{ARCH - Selected readings} (ed R.F. Engle). Oxford:Oxford University Press.

\bibitem[{Fa}{, 2003}]{fa} Fa K. S. (2003) Linear Langevin equation with time-dependent drift
and multiplicative noise term: exact study. {\it Chem. Phys}, \textbf{287}, 1.

\bibitem[{Feder}{, 1988}]{feder} Feder J. (1988) \textit{Fractals}. New York: Plenum.

\bibitem[{Feist}{, 2006}]{feist} Feist A.M., Scholten J.C.M., Palsson B.{\O}., Brockman F.J. and Ideker T. (2006)
Modeling methanogenesis with a genome-scale metabolic reconstruction of M. barkeri. \textit{Mol Syst Biol}, \textbf{2}, 4.

\bibitem[{Finney}{, 1941}]{finney} Finney D. J. (1941) On the distribution of a variate whose logarithm is normally distributed.
\textit{J. Roy. Statist. Soc. B}, \textbf{7}, 155.

\bibitem[{Fr\'{e}chet}{, 1927}]{weibull} Fr\'{e}chet M. (1927) Sur la loi de probabilit\'{e} de l'\'{e}cart maximum.
\textit{Ann. Soc. Pol. Math.}, \textbf{6}, 93.

\bibitem[{Frisch}{, 1997}]{turbulence} Frisch U. (1997) \textit{Turbulence: The Legacy of A. Kolmogorov}. Cambridge:Cambridge University Press.

\bibitem[{Gentle}{, 2004}]{Gentle} Gentle J. E. (2004)
\textit{Random Number Generation and Monte Carlo Methods (Statistics and Computing)}. Berlin: Springer.

\bibitem[{Gibrat}{, 1930}]{gibrat} Gibrat R. (1930) Une loi des r\'{e}partitions \'{e}conomiques.
{\it Bull. Statist. G\'{e}n. Fr.}, \textbf{19}, 469.

\bibitem[{Gradshteyn and Ryzhik}{, 1980}]{gradshteyn} Gradshteyn I. S. and Ryzhik I. M. (1980)
\textit{Table of Integrals, Series, and Products}. New York: Academic Press. \texttt{3.462.1}.
\bibitem[{Green}{, 1948}]{green} Green L. C. (1948) Maximum Uncertainty as a Simple Example of a Non-Distributive Algebra.
{\it Amer. Math Monthly}, \textbf{55}, 363.

\bibitem[{Johnson and Lotz}{, 1970}]{johnsonkotz} Johnson N. L. and Lotz S. (1970)
\textit{Continuous univariate distributions}. New York: John Wiley \& Sons.

\bibitem[{Kapteyn}{, 1903}]{kapteyn} Kapteyn J. C. (1903) \textit{Skew Frequency Curves in
Biology and Statistics}. Groningen: Astronomical Laboratory, Noordhoff.

\bibitem[{Kolmogorov}{, 1941}]{kolmogorov} Kolmogorov A. N. (1941) On the logarithmic
normal distribution law of particles with dimensions of fragmentation. {\it Dok.
Acad. Nauk SSSR}, \textbf{31}, 99.

\bibitem[{Kummel {\it et al.}}{, 2006}]{kummel} K\"{u}mmel A., Panke S. and Heinemann M. (2006)
Systematic assignment of thermodynamic constraints in metabolic network models. \textit{BMC Bioinformatics}, \textbf{7}, 512.

\bibitem[{L\'{e}vy}{, 1954}]{levy} L\'{e}vy P. (1954) \textit{Th\'{e}orie de I'addition des variables al\'{e}atoires}. Paris:Gauthierr-Villards.

\bibitem[{Mandelbrot}{, 1997}]{mandelbrot} Mandelbrot B. B. (1997)
\textit{Fractals and Scaling in Finance}. New York: Springer.

\bibitem[{Nivanen {\it et al.}}{, 2003}]{nivanen} Nivanen L., Le Mehaute A. and Wang Q. A. (2003)
Generalized algebra within a nonextensive statistics. {\it Rep. Math. Phys.},
\textbf{52}, 437.

\bibitem[{Palsson}{, 2006}]{palsson} Palsson B. {\O}. (2006)
\textit{Systems Biology: Properties of reconstructed networks}. Cambridge: Cambridge University Press.

\bibitem[{Petit Lob\~{a}o {\it et al.}}{, 2009}]{thierry} T.C. Petit Lob\~{a}o, P.G.S. Cardoso, S.T.R. Pinho and E.P. Borges (2009)
Some properties of deformed $q$-numbers. {\tt arXiv:0901.4501v1}[math-ph]. Preprint.

\bibitem[{Reed {\it et al.}}{, 2003}]{reed} Reed J.L., Vo T.D., Schilling C.H. and B. {\O}. Palsson (2003) An expanded genome-scale model
of E. coli K-12 (iJR904 GSM/GPR). \textit{Genome Biology}, \textbf{4},R54.1.

\bibitem[{Reed and Palsson}{, 2007}]{reed2} Reed J.L., Palsson B.{\O}. (2007) Genome-Scale in silico models
of E. coli have multiple equivalent phenotypic states: Assessment of correlated reaction subsets that comprise network states. \textit{Genome Res.}, \textbf{14}, 1797.

\bibitem[{Rosin and Rammler}{, 1933}]{weibull1} Rosin P. and Rammler E. (1933) The laws governing the finiteness of Powdered Coal.
\textit{J. Inst. Fuel}, {\bf 7}, 29.

\bibitem[{Stauffer \textit{et al.}}{, 2006}]{murilinho} D. Stauffer, S.M. Moss de Oliveira, P.M.C. de Oliveira and J.M. de S\'{a} Martins (2006)
\textit{Biology, Sociology, Geology by Computational Physicists,}, vol. 1. Amsterdam:Elsevier.

\bibitem[{Tsallis}{, 1988}]{ct-88} Tsallis C. (1988) Possible generalization of
Boltzmann--Gibbs statistics. {\it J. Stat. Phys.}, \textbf{52}, 479.

\bibitem[{Tsallis}{, 1994}]{ct-quimica} Tsallis C. (1994)
What are the numbers that experiments provide? {\it Qu\'{\i}mica Nova}, \textbf{17}, 468.

\bibitem[{Tsallis}{, 2009}]{applications} Tsallis C. (2009) {\it Introduction to
Nonextensive Statistical Mechanics: Approaching a Complex World}. Berlin:Springer.

\bibitem[{Umarov \textit{et al.}}{, 2006}]{clt} Umarov S., Tsallis C., Gell-Mann M. and Steinberg S.,
Symmetric $(q,\alpha )$-Stable Distributions. Part I: First Representation.
\texttt{arXiv:cond-mat/0606038}[cond-mat.stat-mech]. Preprint.
\textit{and} Symmetric $(q,\alpha )$-Stable Distributions. Part II: Second Representation.
\texttt{arXiv:cond-mat/0606040}[cond-mat.stat-mech]. Preprint.

\bibitem[{Umarov and Duarte Queir\'{o}s}{, 2008}]{sudq} Umarov S. and Duarte Queir\'{o}s S. M. (2008) Functional-differential
equations for $F_{q}$-transforms of $q$-Gaussians. \texttt{arXiv:0711.2550}[cond-mat.stat-mech]. Preprint.

\bibitem[{Umarov and Tsallis}{, 2008}]{sabir} Umarov S. and Tsallis C. (2008), On a representation of
the inverse $F_{q}$-transform. {\it Phys. Lett. A} \textbf{372}, 4874.

\bibitem[{von Neumann}{, 1951}]{vonNeuman} von Neumann J. (1951)
Various techniques used in connection with random digits. Monte-Carlo methods. \textit{Nat. Bureau Standards}, \textbf{12}, 36.

\bibitem[{Warren and Jones}{, 2006}]{pbw} Warren P. B. and Jones J. L. (2006)
Duality, Thermodynamics, and the Linear Programming Problem in Constraint-Based Models of Metabolism. \textit{Phys. Rev. Lett.}, \textbf{99}, 108101

\bibitem[{Weibull}{, 1951}]{weibull2} Weibull W. (1951) A statistical distribution function of wide applicability.
\textit{J. Appl. Mech. - Trans. ASME}, \textbf{18}, 293.

\bibitem[{Wolfram}{, 2001}]{paraboliccylinderd} \texttt{ http://functions.wolfram.com/HypergeometricFunctions/ ParabolicCylinderD/}.

\bibitem[{Yuang}{, 1933}]{yuang} Yuan P. T. (1933) On the logarithmic frequency distributions and the semi-logarithmic correlation surface.
\textit{Ann. Math. Statist.}, \textbf{4}, 30.

\end{thebibliography}
\end{document}